\documentclass[%
 aip,
 pop,
 amsmath,amssymb,
 reprint,%
]{revtex4-1}

\usepackage{graphicx}
\usepackage{dcolumn}
\usepackage{bm}

\usepackage[utf8]{inputenc}
\usepackage[T1]{fontenc}
\usepackage{mathptmx}
\usepackage{etoolbox}

\usepackage{amsmath}
\allowdisplaybreaks
\usepackage{amsthm} 
\usepackage[cal=pxtx]{mathalpha}

\newcommand{\vect}[1]{\mathbf{\boldsymbol{#1}}} 
\newcommand{\matr}[1]{\mathbf{#1}} 

\newcommand{\calD}{\mathcal{D}}
\newcommand{\calB}{\mathcal{B}}

\newcommand{\calP}{\mathcal{P}}
\newcommand{\calPm}{\mathcal{P}^{m} }
\newcommand{\calDPm}{\mathcal{DP}^{m}}
\newcommand{\calDXT}{\mathcal{D}\vect{X}_{T}}
\newcommand{\calDXpol}{\calD \vect{X}_{\text{pol}}}

\newcommand{\Bpol}{ \vect{B}_{\text{pol}} }

\newcommand{\rmd}{\mathrm{d}}
\newcommand{\hatb}{\hat{\vect{b}} }

\newcommand{\xcyc}{\vect{x}_{\text{cyc}} }

\newcommand{\vectx}{ \vect{x} }

\usepackage{mathtools} 

\usepackage{scalerel}
\usepackage{stmaryrd}

\def\onedot{$\mathsurround0pt\ldotp$}
\def\cddot{
  \mathbin{\vcenter{\baselineskip.67ex
    \hbox{\onedot}\hbox{\onedot}}%
  }}%
\def\cdddot{
  \mathbin{\vcenter{\baselineskip.67ex
    \hbox{\onedot}\hbox{\onedot}\hbox{\onedot}%
  }}%
}

\usepackage{tikz}
\newcommand{\tikzmark}[1]{\tikz[overlay,remember picture] \node (#1) {};}

\tikzset{down square arrow/.style={to path={-- ++(0,-.25) -| (\tikztotarget)}}}
\tikzset{up square arrow/.style={to path={-- ++(0,+.25) -| (\tikztotarget)}}}
\tikzset{
    up square arrow cdot/.style={
        to path={-- ++(0,0.25) -| (\tikztotarget) node[pos=0.25, fill=white] {$\cdot$}}
    },
    down square arrow cdot/.style={
        to path={-- ++(0,-0.45) -| (\tikztotarget) node[pos=0.25, fill=white] {$\cdot$}}
    }
}

\usepackage{cancel}

\makeatletter
\def\@email#1#2{%
 \endgroup
 \patchcmd{\titleblock@produce}
  {\frontmatter@RRAPformat}
  {\frontmatter@RRAPformat{\produce@RRAP{*#1\href{mailto:#2}{#2}}}\frontmatter@RRAPformat}
  {}{}
}%
\makeatother
\begin{document}

\preprint{AIP/123-QED}

\title[Shifts of orbits and periodic orbits]{On the shifts of orbits and periodic orbits under perturbation and \\the change of Poincaré map Jacobian of periodic orbits}
\author{Wenyin Wei}
     \affiliation{%
    Institute of Plasma Physics, Hefei Institutes of Physical Science, Chinese Academy of Sciences, Hefei 230031, People's Republic of China
    }
    \affiliation{
    University of Science and Technology of China, Hefei 230026, People's Republic of China
    }
    \affiliation{
    Forschungszentrum J\"{u}lich GmbH, Institute of Fusion Energy and Nuclear Waste Management – Plasma Physics, 52425 J\"{u}lich, Germany
    }
\author{Alexander Knieps}
    \affiliation{
    Forschungszentrum J\"{u}lich GmbH, Institute of Fusion Energy and Nuclear Waste Management – Plasma Physics, 52425 J\"{u}lich, Germany
    }    
\author{Yunfeng Liang*$^{,}$}
    \email{y.liang@fz-juelich.de}
    \affiliation{%
    Institute of Plasma Physics, Hefei Institutes of Physical Science, Chinese Academy of Sciences, Hefei 230031, People's Republic of China
    }
    \affiliation{
    Forschungszentrum J\"{u}lich GmbH, Institute of Fusion Energy and Nuclear Waste Management – Plasma Physics, 52425 J\"{u}lich, Germany
    }

\date{\today}

\begin{abstract}
Periodic orbits and cycles, respectively, play a significant role in discrete- and continuous-time dynamical systems~(\textit{i.e.} maps and flows). To succinctly describe their shifts when the system is applied perturbation, the notions of functional and functional derivative are borrowed from functional analysis to consider the whole system as an argument of the geometric representation of the periodic orbit or cycle. The shifts of an orbit/trajectory and periodic orbit/cycle are analyzed and concluded as formulae for maps/flows, respectively. The theory shall be beneficial for analyzing sensitivity to perturbations, and optimizing and controlling various systems. 
\end{abstract}
\keywords{periodic orbit, magnetic topology, functional perturbation theory, chaotic map, functional derivative}
\maketitle


\small

\section{Introduction}
This paper develops a theory on the \textit{shifts of orbits and periodic orbits under perturbation}, motivated by the practical need for three-dimensional (3D, \textit{i.e.} non-axisymmetric) magnetic field design optimization and control in magnetically confined fusion (MCF) machines. The theory even applies to arbitrary finite-dimensional dynamical systems (DSs) generated by ordinary differential equations (ODEs). Generalization to infinite-dimensional DSs~\cite{rudnicki2004} generated by partial differential equations~(PDE) is possible but requires more effort beyond this Letter, which can offer a convenient computational method to estimate the sensitivity of the system to the initial conditions.

In the MCF community, the past practice~\cite{yunfeng2007, yunfeng2010, yunfeng2013, frerichs2020, roeder2003, evans2002, evans2005}, to mitigate or suppress Edge Localized Modes~(ELMs) in tokamak plasmas involves introducing external Resonant Magnetic Perturbations~(RMPs) to let island chains grow and corrode each other, then a chaotic field is created, also referred to as a stochastic field \cite{abdullaev2015, abdullaev2016magnetic, howard1984} owing to an incomplete understanding of the long-term ergodic behaviour of field line tracing. It was expected that the chaotic field between island chains would reduce the transient heat flux released by ELMs by enhancing radial transport at the plasma edge. However, the theory underlying this practice has shortcomings due to relying on Fourier analysis of the radial magnetic perturbation $B^r = \vect{B}\cdot\nabla r$ in the flux coordinates $(r, \theta, \varphi)$ or its equivalents~\cite{balescu1998tokamap, balescu1998} that require the integrability of the perturbed system. Although the term \emph{stochastic field} is still widely used in the MCF community and Hamiltonian system-relevant research~\cite{unterberg2010, mackay1984}, we advocate using the term \emph{chaotic field} to avoid confusion, as researchers sometimes use chaotic and stochastic interchangeably unconsciously or consciously. There exists a mathematical branch of stochastic process analysis, which is clarified to be irrelevant to the chaos discussed here. The chaotic field expanded by the stable and unstable manifolds of the outermost X-cycle(s) is probably the largest one in MCF machines~\cite{frerichs2015}, prompting a return to the real-world cylindrical coordinates to avoid errors induced by transforming.


The chaotic field is not unique in toroidal magnetic fields but is common in autonomous continuous-time DSs of three dimensions or more, not necessarily divergence-free. Chaos in 2D continuous-time DSs is forbidden by Poincar\'e-Bendixson theorem. For discrete-time DSs, one dimension is sufficient to allow chaos to exist. 

Functional, a function of function, is not a novel concept introduced to physicists from mathematics. In quantum physics, it is a common practice to write various forms of energy as functions of the electronic wave function $\psi$, which has established a systematical framework of density functional theory~(DFT). Considering the whole system as an argument of orbits and periodic orbits, the \textit{functional perturbation theory}~(FPT) developed in this Letter significantly enhances one's ability to analyze the change of these well-defined objects under perturbation, which can be a bit more complicated than DFT due to the more complex composite relationships between functions.

Apart from MCF, this theory has a broad application range across numerous research domains concerned with the behaviour of a dynamical system~\cite{harsoula2018, meiss2015, zotos2015, henok2022, falessi2015, digiannatale2018, pegoraro2019, serra2023, memarian2024} or the island-around-island hierarchy in the chaotic field~\cite{meiss1986, alus2014, abdullaev2016magnetic}. With aircraft orbit optimization in aerospace dynamics as an example\cite{hapala2016, trelat2012, runqi2019}, the propulsion force contributed by the engine can be treated as a kind of perturbation to the system. Then, when and in which direction to expel the exhaust gases in pursuit of the shortest time, the highest fuel efficiency, the heaviest load or the best trade-off among them become a problem that can be solved by leveraging this theory.

\section{Deduction and demonstration}
For an $N$-dimensional~($N$-D) dynamical system in Cartesian coordinates and a 3D toroidal vector field in cylindrical coordinates as a special case of the former, as shown below,
\begin{subequations}
\begin{align}
  \dot{  \vect{x}  } &= \vect{B}(\vect{x}) , \\
  \dot{  \vect{x} }_{\text{pol}} &=
  \frac{ R \vect{B}_{\text{pol}} }{ B_\phi } (\vect{x}_{\text{pol}}, \phi) ,
\end{align}
\label{eq:fundamental_tracing}
\end{subequations}
write a trajectory as $\vect{X}(\vect{x}_0, t)$ and $\vect{X}(\vect{x}_{0,\text{pol}}, \phi_s, \phi_e )$, where $\vect{x}_{0}$ and $\vect{x}_{0,\text{pol}}$ are initial conditions, $\vect{B}_{\text{pol}}$ and $B_\phi$ are poloidal and toroidal components of the field in standard cylindrical coordinates, $\phi$ the azimuthal angle, $\phi_s$ and $\phi_e$ the starting and ending angles.

The theory is inspired by chaotic field-relevant research in the MCF community, so the deduction and results for these two forms are displayed together to facilitate utilizing them for the audience interested in either one. One thing worth emphasizing is that the theory in this Letter does not require the field $\vect{B}$ to be divergence-free, so it can be applied to general vector fields, \textit{e.g.} fluid velocity and current density fields. With derivatives in the initial conditions, equations \eqref{eq:fundamental_tracing} become
\begin{subequations}
\begin{align}
  \frac{\partial }{\partial t } \mathcal{D}\vect{X}
  & = 
    \nabla \vect{B} \cdot \mathcal{D} \vect{X}
    \label{eq:DX_timeforward} \\
  \frac{\partial }{\partial \phi_e } \calDXpol
  &=
  \underbrace{
    \frac{\partial (R\vect{B}_{\text{pol}}/B_{\phi ) } }{\partial (R,Z)}
    }_{\mathrlap{\smash{ \text{abbr. as } \matr{A}=\matr{A}(R,Z,\phi) } }} 
    \cdot \calDXpol  
   \label{eq:DXpol_timeforward}
\end{align}
\end{subequations}
where 
$\mathcal{D}$ means the partial derivative \textit{w.r.t.} $\vect{x}_0$ or $\vect{x}_{0,\text{pol}}$. Note that $\matr{A} = \matr{A}(R, Z,\phi)$ is a function of spatial coordinates directly, $\vect{X}_\text{pol}(\vect{x}_{0,\text{pol}}, \phi_s, \phi_e)$ is a function of the initial point (not \emph{in-situ}). To emphasize the difference, two derivative operators $\partial / \partial (R,Z)$ and $\mathcal{D}$ are used.

\begin{figure}
\includegraphics[width=0.65\linewidth]{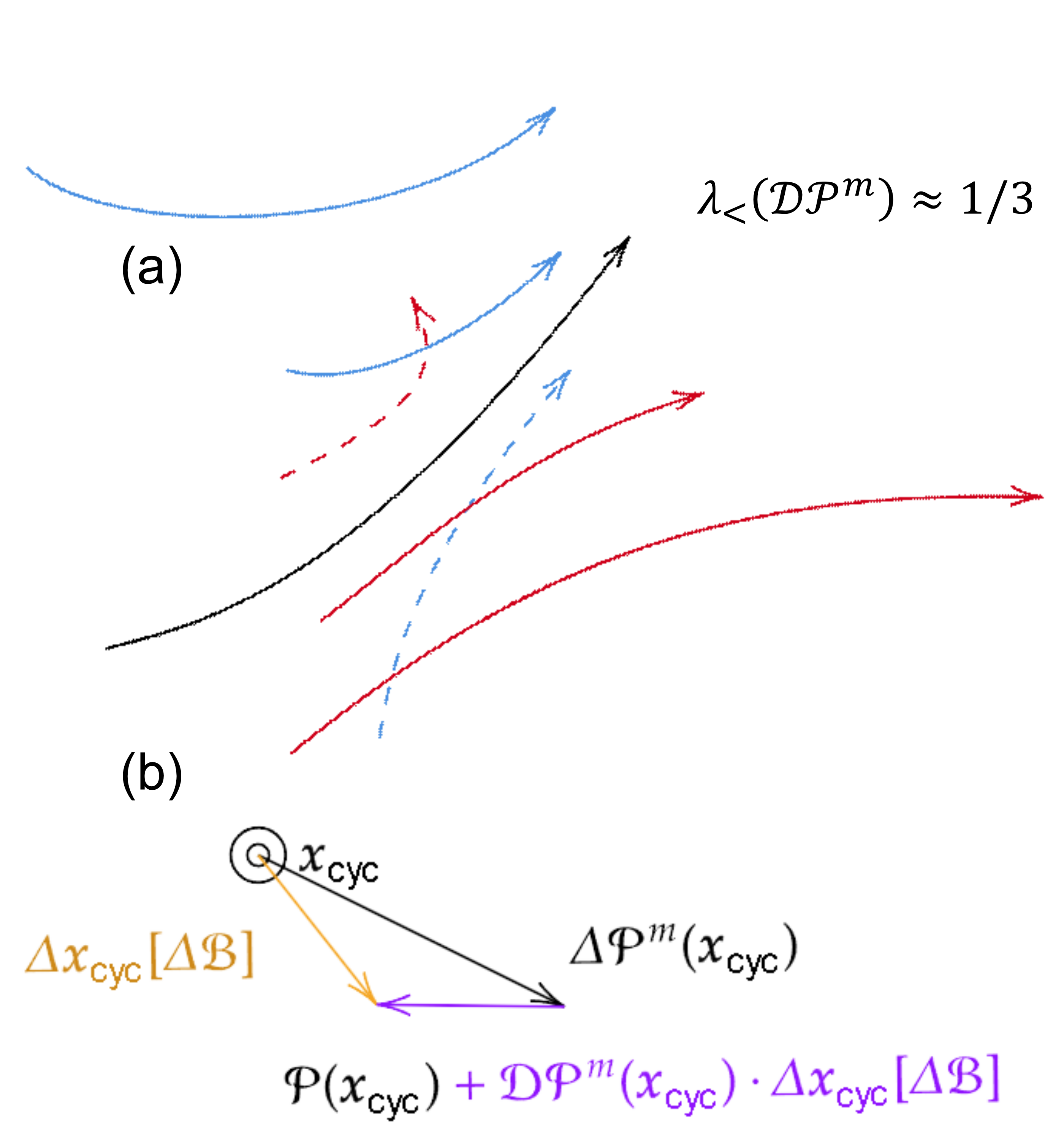}
\caption{\label{fig:DPm_relevant_cartoon} (a) Cartoon to show the geometry meaning of $\calDPm$ eigenvalues, (b) Cartoon to show how to calculate the shift of a fixed point for a map under perturbation. }
\end{figure}

The term X-cycle is a figurative alias of \textit{hyperbolic cycle} (see Fig.~\ref{fig:DPm_relevant_cartoon}(a)). A cycle is \textit{hyperbolic} if 
\begin{align}
 &|\lambda_{i}|\neq 1 \text{ for all eigenvalues of } \calDPm , \nonumber \\
 \text{or equiv. }  &|\lambda_{i}|\neq 1 \text{ for all eigenvalues of } \calDXT \text{ except} \nonumber
\end{align}
the one corresponding to the eigenvector $\vect{v}_{i}=\hatb$, where $\hatb$ is the local field normalized. An O-cycle can be similarly defined in 3-D continuous-time systems with both $\calDPm$ eigenvalues $\lambda_i \in \mathbb{S}\subset \mathbb{C}$ but $\neq \pm 1$. $\calP$ denotes the Poincar\'e map. For a 3D toroidal vector field, the Poincar\'e section is chosen to be, in this Letter, an $R$-$Z$ cross-section, while $m$ is the toroidal turn number of the cycle. For a general $N$-D system, $\calPm$ can be defined as a whole with no individual meaning of $m$ to express the returning map for quasi-periodic cases.
The evolution of $\calDPm$ and $\calDXT$ along a cycle has been revealed by two formulae given in \cite{wei2023} as shown below: 
\begin{subequations}
\begin{align}
  \frac{\mathrm{d} }{\mathrm{d} t } \mathcal{D}\vect{X}_T &= \left[
    \nabla \vect{B}, \mathcal{D} \vect{X}_T
  \right]
  \label{eq:DXT_evolution} \\
\intertext{where $\mathcal{D}\vect{X}_T$ is short for the full-period Jacobian of the cycle, that is $\mathcal{D}\vect{X}( \vect{x}_0, T)|_{\vect{x}_0 = \vect{X}(\vect{x}_{\text{fix}},t) }$, $\vect{x}_{\text{fix}}$ is a fixed point on the cycle, and}
  \frac{\mathrm{d} }{\mathrm{d} \phi }\calDPm 
  &= \left[
        \frac{\partial (R\vect{B}_{\text{pol}}/B_{\phi ) } }{\partial (R,Z)}, 
    \calDPm
  \right] = \left[ \matr{A}, \calDPm \right]. 
  \label{eq:DPm_evolution}
\end{align}
\end{subequations}


Partial and total functional derivatives are denoted by $\delta / \delta \mathcal{B} $ and $\rmd / \rmd \mathcal{B} $, respectively, but they are inconvenient to use due to their infinite-dimensional nature. Hence, they are often accompanied by a given perturbation $\Delta\mathcal{B}$ to be directional derivatives $\Delta \mathcal{B} \cdot \delta / \delta \mathcal{B} $ and $\Delta \mathcal{B} \cdot \rmd / \rmd \mathcal{B} $. In this Letter, when a function symbol is written in calligraphic font, it means that it is not to be evaluated at a specific point but considered as a standalone object. The Poincar\'e map $\mathcal{P}$ is a special case also written in calligraphic font but for not to be misunderstood as a vector filed.
An alternative notation of directional partial derivative equivalent to $\prod_k (\Delta \calB_k \cdot \delta / \delta\calB) \vect{X}$ is $\delta^k \vect{X}[\calB; \Delta \calB_1, \dots, \Delta \calB_k ]$ ~(see \cite{frigyik2008}), of which the argument list can be omitted for brevity. 

If the system itself is evaluated at a specific point and applied the functional derivative with always the same $\Delta \calB$, the result is denoted as $\delta^k \vect{B}[\calB;\Delta \calB, \dots, \Delta\calB](\vect{x})$. Collecting these derivatives across all $\vect{x}$ in the domain of $\calB$ allows one to construct a new vector function, denoted by $\delta^k \calB$. Note that the denotation $\delta^k \calB$ can be conveniently expressed without the function argument list $[\calB; \Delta\calB, \dots, \Delta\calB]$ when the context is clear. A perturbation $\Delta \calB$ then can be decomposed into different orders as follows,
\begin{align*}
\calB + \Delta \calB = \calB + \delta\calB + \frac{1}{2!}\delta^2 \calB
+ \frac{1}{3!}\delta^3 \calB + \cdots.
\end{align*}
Conversely, one can also define a perturbation $\Delta\calB$ by specifying a sequence of $\delta^k \calB$. 
For example, 
\begin{align*}
\vect{J} \times \vect{B}  =& \nabla p, \\
\delta\vect{J} \times \vect{B} + \vect{J} \times \delta\vect{B}  =& \delta \nabla  p, \\
\delta^2\vect{J} \times \vect{B} + 2 \delta\vect{J} \times \delta\vect{B} + \vect{J} \times \delta^2 \vect{B} =& \delta^2 \nabla  p, \\ 
\quad\cdots
\end{align*}
Then one can write $(\vect{J}+\Delta\vect{J}) \times (\vect{B}+\Delta\vect{B}) = \nabla p + \Delta\nabla p
$ with $\Delta\vect{J}$, $\Delta\vect{B}$ and $\Delta\nabla p$ defined by their corresponding series expansions. 

With the whole system $\calB$ considered as an argument of $\vect{X}$, the equation~(\ref{eq:DX_timeforward}) becomes
\begin{equation}
\frac{\partial }{\partial t} \vect{X}\left[ \calB \right] (\vect{x}_0, t) 
= \vect{B} (\vect{X}) 
= \vect{B} \left[ \calB \right] (\vect{X} \left[\calB \right] (\vect{x}_0, t) ) .
\label{eq:X_timeforward_variantB}
\end{equation}
Impose directional functional derivatives on both sides, then
\begin{align}
\frac{\partial }{\partial t}  \delta \vect{X}\left[ \calB; \Delta \calB \right] (\vect{x}_0, t) 
&= \delta \vect{B} \left[ \calB ; \Delta \calB \right] \left( \vect{X} [ \calB ] (\vect{x}_0, t) \right)  \nonumber \\
&+
\delta \vect{X} \left[\calB ; \Delta \calB \right] (\vect{x}_0, t) 
\cdot 
\nabla \vect{B} 
\label{eq:deltaX_timeforward} 
\end{align}
Since the the perturbation $\Delta \calB$ is usually a given one rather than a series of different perturbations $(\Delta \calB_1, \Delta \calB_2, \dots)$,  $\delta \vect{X} \left[\calB; \Delta\calB \right](t)$ can be concisely represented as $\delta\vect{X}$, similarly $\delta^2 \vect{X} \left[\calB; \Delta\calB, \Delta \calB \right](t)$ as $\delta^2 \vect{X}$. With arguments omitted, the equation~(\ref{eq:deltaX_timeforward}) becomes 
\begin{equation}
\frac{\partial }{\partial t}  \delta \vect{X}
= \delta \vect{B} 
+  \delta \vect{X} \cdot \nabla \vect{B} 
\label{eq:progress_deltaX_cartesian}
\end{equation}

Similarly, for the second and third variations,
\begin{equation}
\frac{\partial }{\partial t}  \delta^2 \vect{X}
= \delta^2 \vect{B} + 2  \left( \delta \vect{X} \cdot \nabla \right) \delta \vect{B} 
+ \delta \vect{X}  \delta \vect{X} \cddot \nabla^2 \vect{B}
+ \delta^2 \vect{X}  \cdot \nabla \vect{B} .
\label{eq:delta2X_timeforward}
\end{equation}
where the semicolon is defined as $\vect{ab} \cddot \vect{cd}= (\vect{a}\cdot\vect{c}) (\vect{b}\cdot\vect{d})$,
\begin{align}  
\frac{\partial }{\partial t}  \delta^3 \vect{X} 
& = \delta^3 \vect{B} + 3 (\delta \vect{X}\cdot \nabla ) \delta^2  \vect{B}
\label{eq:delta3X_timeforward} \\
& + 3 (\delta^2\vect{X}\cdot\nabla) \delta \vect{B}   
+ 3 (\delta\vect{X}\delta\vect{X}\cddot\nabla^2) \delta \vect{B} 
\nonumber\\
&+  3 \delta\vect{X}\delta^2\vect{X}\cddot\nabla^2 \vect{B} 
+ \delta\vect{X}\delta\vect{X}\delta\vect{X} \cdddot \nabla^3 \vect{B} 
+ \delta^3 \vect{X}\cdot \nabla \vect{B} \nonumber
\end{align}

The deduction can be continued to higher orders. The shift of a trajectory under perturbation to the whole system can be approximated by the following Taylor expansion based on these variations,
\begin{align}
\vect{X} \left[\calB+\Delta\calB \right](\vect{x}_0, t) 
&=
\vect{X}\left[\calB  \right](\vect{x}_0, t)  
+ \delta \vect{X}    \nonumber \\
+ \delta^2 \vect{X}/2! + \cdots &+ \delta^k \vect{X}/k! + \mathcal{O}(\|\Delta \calB \|^{k+1} ) .
\label{eq:DeltaX_Taylor_expansion}
\end{align}

\begin{figure}
    \includegraphics[width=0.8\linewidth]{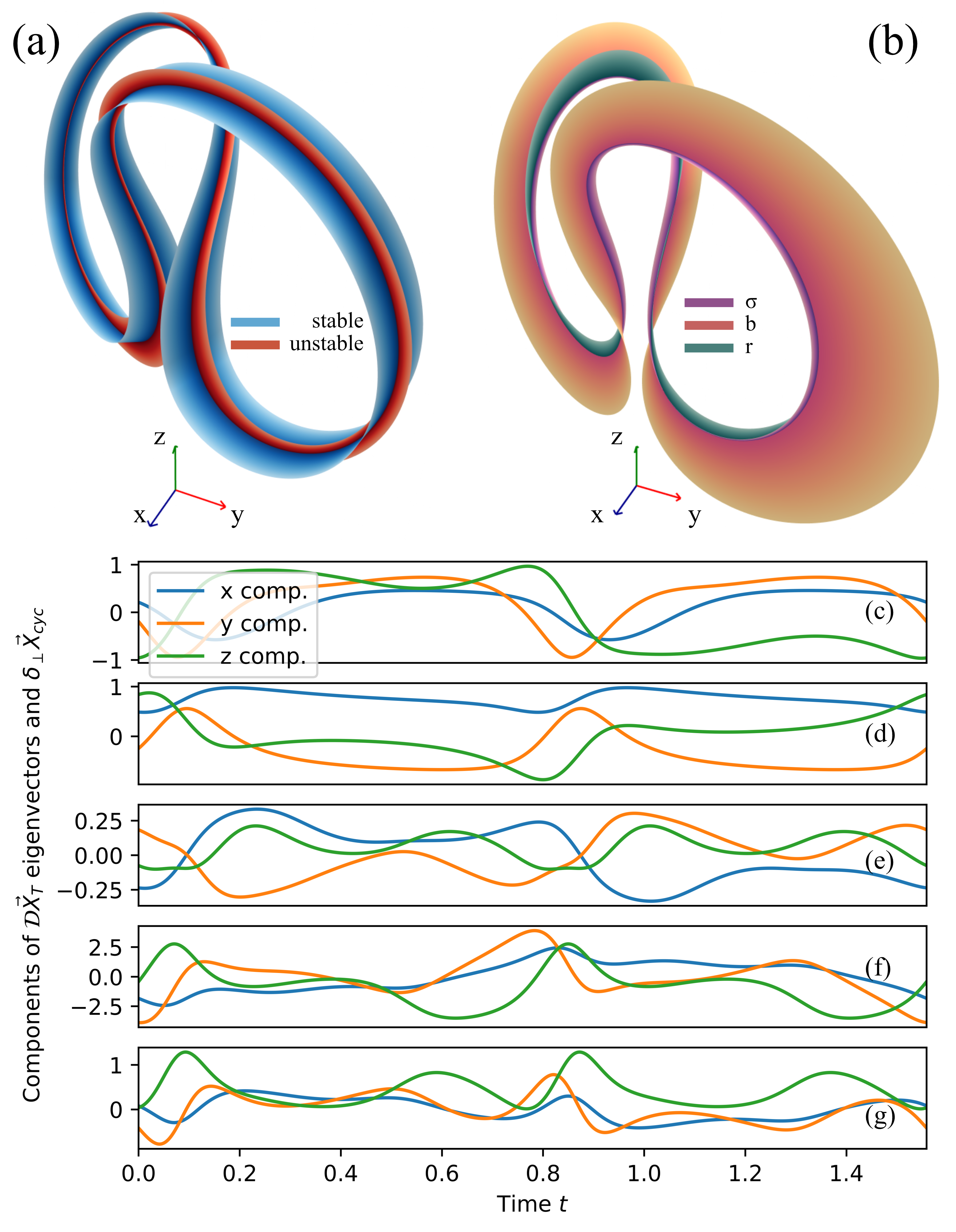}
    \caption{\label{fig:Lorenz_ABcycle_fig_DXT_delta_perp_X_cyc} A typical cycle of Lorenz attractor with period $T=1.5586522107162$, $\vect{x}_0 = [-13.763610682134,~ -19.578751942452,~ 27]$~(chosen from Table 1 in \cite{viswanath2003}). (a and c~-~d) $\calDXT$ eigenvectors (except the one $\vect{v}_{i}=\hatb$) and their components. (b and e~-~g) $\delta_{\perp}\vect{x}_{\text{cyc}} [\mathcal{B}; \Delta \mathcal{B}], ~ \Delta\mathcal{B}\in [\mathcal{B}_{\sigma}, \mathcal{B}_{b}, \mathcal{B}_{r}] $ and their components.}
    \includegraphics[width=1.0\linewidth]{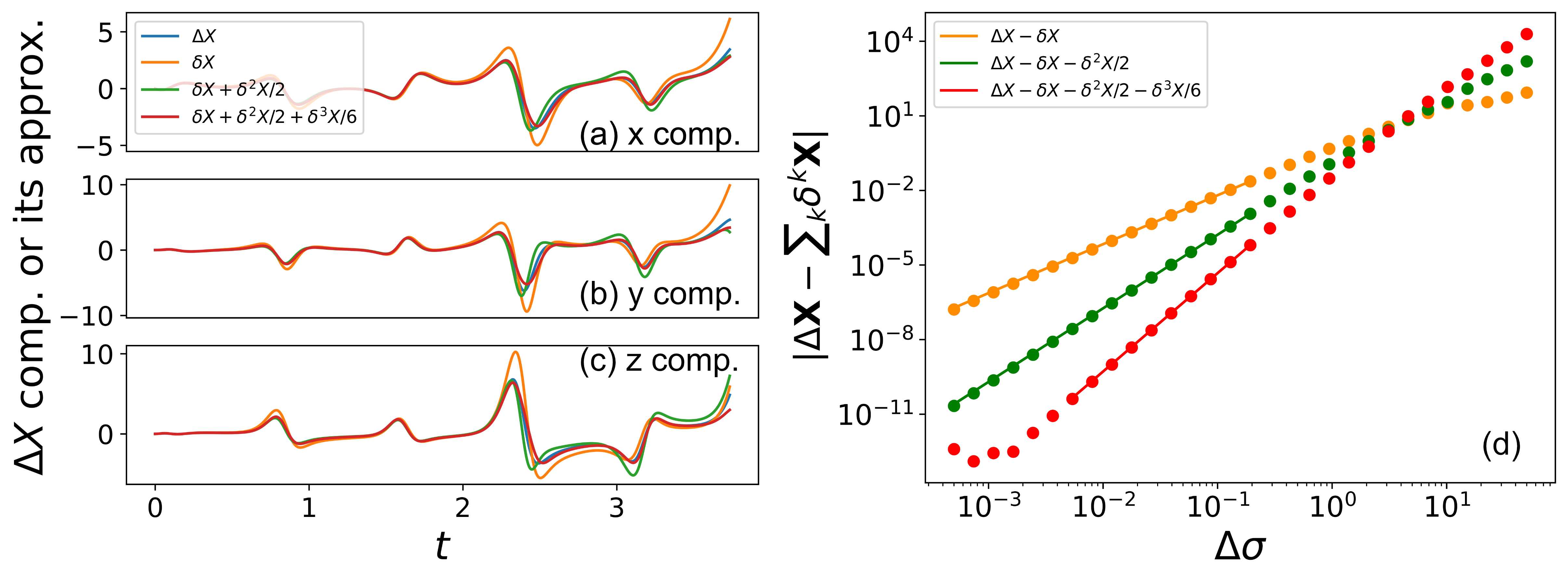}
    \caption{\label{fig:Lorenz_deltaX_approx} (a) The components of $\Delta\vect{X}$ with $\Delta\sigma=1.0$ and its approximations. (b) The approximation errors against $\Delta \sigma$ at $t=1.0$, .}
\end{figure}
Take Lorenz attractor (Fig.~\ref{fig:Lorenz_ABcycle_fig_DXT_delta_perp_X_cyc} and \ref{fig:Lorenz_deltaX_approx}) as an example to show the power law of error. Set $\sigma=10$, $r=28$, and $b=8/2$, a standard choice of the parameter set. Change $\sigma$ to perturb the system, as shown below
\begin{align*}
    \vect{B} = 
\begin{bmatrix}
\sigma (-x+y)\\ 
-xz +rx -y \\ 
xy - bz
\end{bmatrix},
\quad
\Delta_\sigma \vect{B} = 
\begin{bmatrix}
-x+y\\ 
0 \\ 
0
\end{bmatrix}.
\end{align*}
Then one shall expect the approximation error $|\Delta\vect{X}-\sum _{k} ~ \delta^k \vect{X}|$ to be an order of $\mathcal{O}(\Delta \sigma^{k+1})$, where $\Delta\vect{X}:=\vect{X}[\mathcal{B}+\Delta\mathcal{B}] - \vect{X}[\mathcal{\mathcal{B}}] $. Then 
\begin{align} 
\ln  \text{err} & \propto  k_{\text{err}} \ln \Delta \sigma, \quad k_{\text{err}} \geq k+1,
\end{align}
which is verified in Fig.~\ref{fig:Lorenz_deltaX_approx}(b) by fitting the scatters. The fitted slopes for the first three order approximations, as shown by straight lines in Fig.~\ref{fig:Lorenz_deltaX_approx}(b), equal 1.96772, 2.96082, 3.95880, respectively. The fitted slopes coincide well with $k+1$ as expected.

To migrate Eqs.~(\ref{eq:X_timeforward_variantB}-\ref{eq:delta3X_timeforward}) from flows to maps, \textit{e.g.} a general map $\calP: \mathbb{R}^N \to \mathbb{R}^{N}$, one simply needs to replace $\partial \vect{X}/\partial t$ on the left hand side~(LHS) with $\calP^{k+1}(\vectx)$, and $\vect{X}$ on the right hand side~(RHS) with $\calP^{k}(\vectx)$. For example, Eq.~\eqref{eq:X_timeforward_variantB} becomes
\begin{equation}
\calP^{k+1} \left[\calP \right] (\vectx) 
= \calP [\calP] (\calP^k \left[\calP \right] (\vectx) ),   
\end{equation}

Equations~(\ref{eq:X_timeforward_variantB}-\ref{eq:DeltaX_Taylor_expansion}) do not need to make substantial change if the system is non-autonomous, \textit{i.e.} the time is an explicit variable of the field $\vect{B}=\vect{B}(t,\vect{x})$.

Thereafter, one can calculate the shifts of a periodic orbit and a cycle \textit{resp.} based on the first variations $\delta \mathcal{P}^m$ and $\delta \vect{X}_T$ with the aid of $\calDPm$, as explained in Fig.~\ref{fig:DPm_relevant_cartoon}(b). Note that a periodic orbit and cycle are defined to be those the initial point of which is tied to the ending point. The idea to calculate the shift is $\vect{x}_{\text{cyc}}$ must have such a shift $\Delta \vect{x}_{\text{cyc}}$ that the ending point $\mathcal{P}(\vect{x}_{\text{cyc}} )$ after perturbation keeps matching the starting point, \textit{i.e.}
\begin{align}
\Delta \vect{x}_\text{cyc} 
&= \Delta\mathcal{P}^m (\vect{x}_\text{cyc}) 
+ \calDPm(\vect{x}_\text{cyc})\cdot \Delta \vect{x}_\text{cyc} + \cdots ,
\label{eq:xcyc_shift_how}\\
\Delta \vect{x}_\text{cyc} 
&= 
\left[
    \calDPm(\vect{x}_\text{cyc}) - \matr{I}
\right]^{-1} \cdot
\left[
 - \Delta \mathcal{P}^{m}(\vect{x}_\text{cyc})  
\right] \nonumber\\
&=\left[
    \calDPm(\vect{x}_\text{cyc}) - \matr{I}
\right]^{-1} \cdot \nonumber \\
\Big( & 
    -\delta \mathcal{P}^m [\calB;\Delta \calB](\vect{x}_\text{cyc}) 
    + \mathcal{O}(\|\Delta \calB\|^2) 
\Big) + \cdots, \\
\intertext{$\delta\mathcal{P}^m$ on RHS contributes the first order variation, so} 
\delta\vect{x}_{\text{cyc}} 
&= - \left[
    \calDPm - \matr{I}
\right]^{-1} \cdot \delta \mathcal{P}^m ,
\label{eq:DeltaXcyc}
\end{align}
which applies to periodic orbits and cycles whose full-period Jacobians $\calDPm$ do not have eigenvalues equal to one. For those on center manifolds, this condition is broken, requiring a different analysis. Deduction for higher-order $\delta^k\xcyc$ is put in Appendix \ref{ap:high_order_shifts}. 



\begin{figure*}[htbp]
    \includegraphics[width=0.75\linewidth]{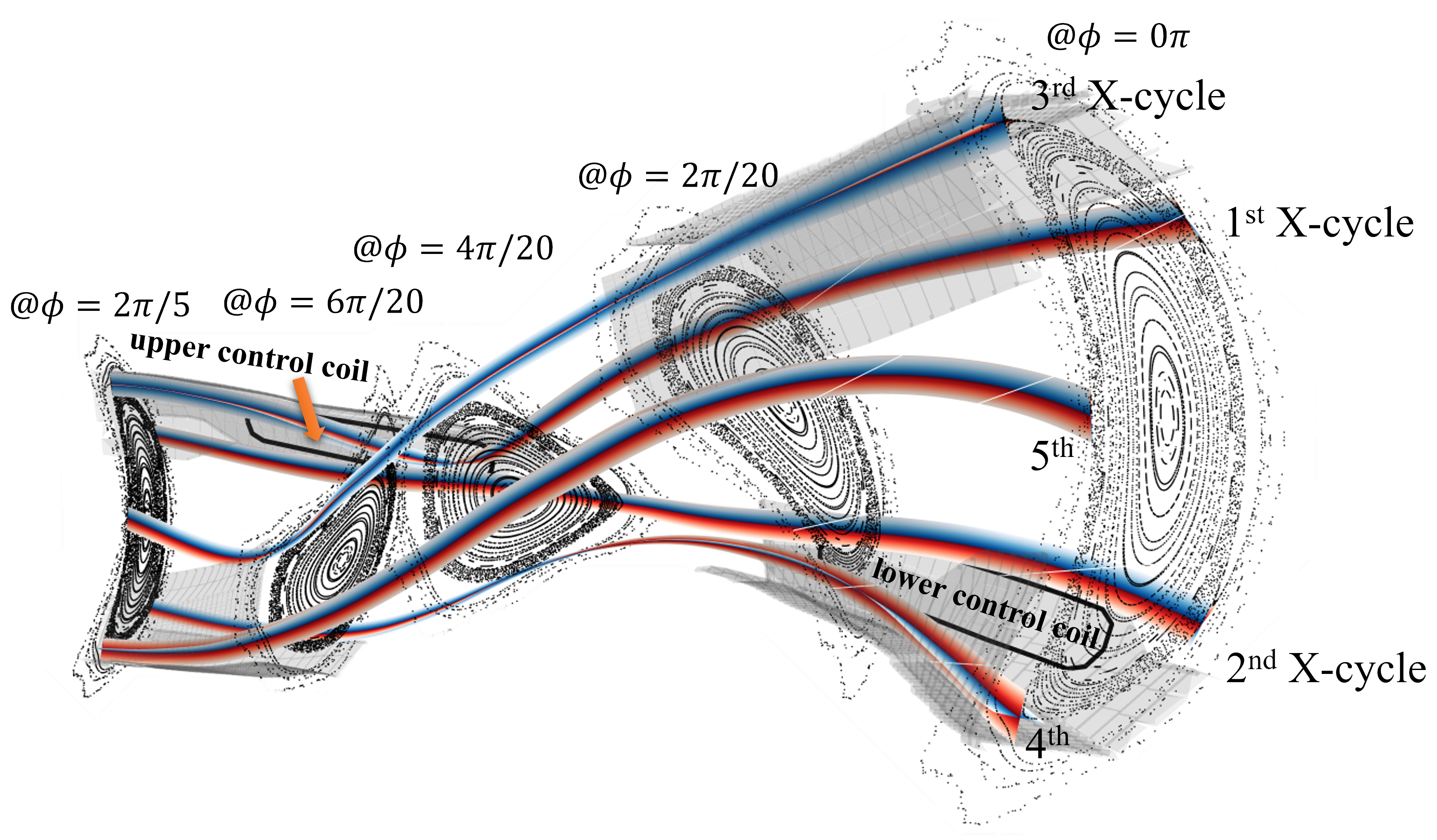}
    \caption{\label{fig:W7X_Poincare_3D} Five typical Poincar\'e plots for Wendelstein 7-X standard configuration. Red and blue ribbons \textit{resp.} indicate the directions of $\mathcal{DP}^1$~($m=1$) unstable and stable eigenvectors corresponding to $\lambda_{i}=1.94965374$ and $0.51291252$.}
\end{figure*}

\begin{figure}[htbp]
    \includegraphics[width=\linewidth]{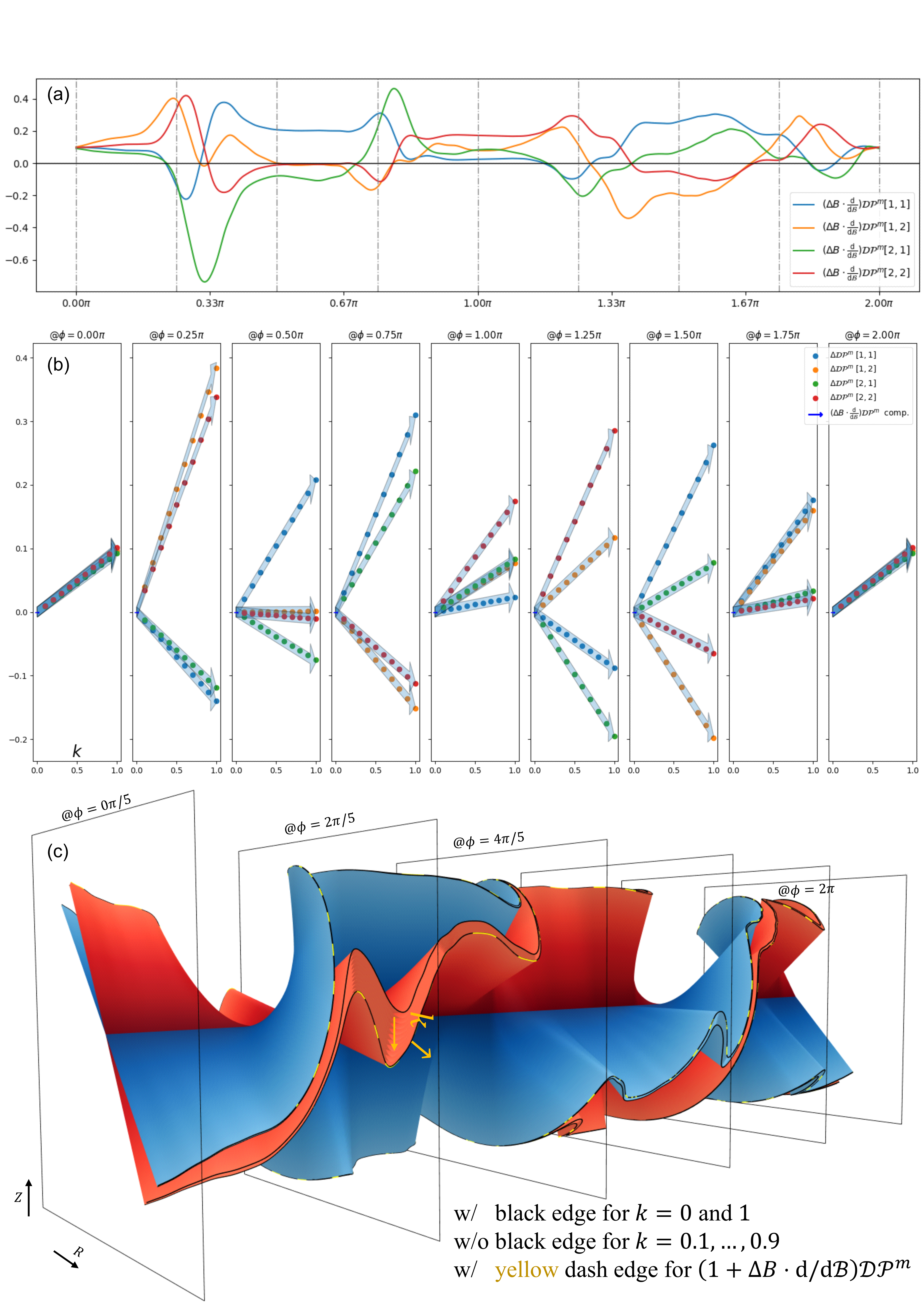}
    \caption{\label{fig:DPm_change_Xcyc_of_W7X} The change of $\calDPm$ of the first X-cycle at the edge $\iota = n/m= 5/5$ island chain of Wendelstein 7-X standard configuration. $(\Delta\calB\cdot \rmd / \rmd \calB) \calDPm$ components are shown as curves in (a) and as transparent blue arrows in (b). Scatter points in (b) show the reference values for the changes of $\calDPm$ components computed for each $k$ in a brute-force way. (c) $\calDPm$ eigenvectors against $\phi$. }
\end{figure}
To apply Eq.~(\ref{eq:DeltaXcyc}) to an $N$-D flow, one needs to define a local $(N-1)$-D Poincar\'e plane for every point on the concerned cycle to define Poincar\'e map. A natural choice is the local plane perpendicular to the cycle. To avoid constructing one more frame of coordinates on the plane, a suggested way is to utilize the relationship between $\calDXT$ and the local Poincar\'e map Jacobian $\calDPm$. $\calDPm$ has the same eigenvalues $\lambda_i$ as $\mathcal{D}\vect{X}_{T}$ but with eigenvectors projected to the Poincar\'e plane $\vect{v}_{i\perp}=\vect{v}_i- \hatb\hatb^T \cdot \vect{v}_i$. Decompose $\Delta \vect{x}_\text{cyc}$ to be the sum of eigenvectors as $\Delta \vect{x}_\text{cyc}\approx \delta_{\perp} \vect{x}_{\text{cyc}} = \sum c_i \vect{v}_{i\perp}$. Then one can replace $\calDPm(\vect{x}_\text{cyc})\cdot \Delta \vect{x}_\text{cyc}$ in Eq.~(\ref{eq:xcyc_shift_how}) with $\sum c_i \lambda_i \vect{v}_{i\perp}$ and $\Delta\mathcal{P}^m(\vect{x}_{\text{cyc}})  $ with its first variation $\delta_{\perp} \vect{X}_T:= \delta\vect{X}_T - \hatb\hatb^T \cdot \delta\vect{X}_T$. Then one can solve for the coefficients $c_i$ to determine $\delta_{\perp}  \vect{x}_{\text{cyc}}$. The other way is to reuse the definition, that is keeping the starting and ending points tied to solve for the shift, but this time take into account projecting $\delta \vect{X}_T$ and $\calDXT\cdot \delta_{\perp} \vect{x}_{\text{cyc}}$ to the perpendicular Poincar\'e plane as below,
\begin{equation}
    \delta_{\perp} \vect{x}_{\text{cyc}} 
    = (\matr{I}-\hat{\vect{b}}\hat{\vect{b}}^T)\cdot \delta \vect{X}_T 
    + (\matr{I}-\hat{\vect{b}}\hat{\vect{b}}^T)\cdot \calDXT \cdot \delta_\perp \vect{x}_{\text{cyc}} ,
    \label{eq:delta_perp_xcyc_by_definition}
\end{equation}
of which the LHS is the starting point shift, while the RHS is the ending point shift that consists of two terms: the first one is contributed by the perturbation and the second one is propagated from the starting point shift via $\calDXT$. Reorganize Eq~\eqref{eq:delta_perp_xcyc_by_definition} to be
\begin{equation}
\delta_{\perp} \vect{x}_{\text{cyc}}  
= 
\left[
    \matr{I} - (\matr{I} - \hat{\vect{b}}\hat{\vect{b}}^T) \cdot \calDXT
\right]^{-1} \cdot
(\matr{I}-\hat{\vect{b}}\hat{\vect{b}}^T)\cdot \delta \vect{X}_T .
\end{equation}

To avoid repetitive calculation of $\delta \vect{X}_T$ and $\delta \mathcal{P}^m$ for every point on a cycle, their evolution along the cycle is concluded as the following formulae,
\begin{subequations}
\begin{align}
    \frac{\rmd}{\rmd t}\delta \vect{X}_T 
    & = \nabla\vect{B} \cdot \delta \vect{X}_T - (\calDXT - \matr{I}) \cdot \delta \vect{B}, 
    \label{eq:deltaXT_evolve_Cartesian}
    \\
    \frac{\rmd}{\rmd \phi}\delta \mathcal{P} ^m
    & = \frac{\partial (R\Bpol / B_\phi)}{\partial (R,Z)} \cdot \delta \mathcal{P} ^m - (\calDPm - \matr{I}) 
    \cdot \delta  \frac{ R\vect{B}_{\text{pol}} }{ B_\phi },  
    \label{eq:deltaPm_evolve_cylindrical}
\end{align}      
\end{subequations}
where $\delta (R\vect{B}_{\text{pol}} / B_\phi )$ is short for $(\Delta \calB \cdot \delta / \delta \calB) \left(
R \vect{B}_{\text{pol}} / B_\phi
\right) $, equal to $\left(
        \dfrac{R~\delta \Bpol }{B_\phi} 
        - \dfrac{ R \Bpol }{ B_\phi^2 }  \cdot \delta B_\phi
    \right) $ by the product rule of differentiation.

Similarly, the calculation of $\delta \vect{x}_{\text{cyc}}$ given by Eq.~(\ref{eq:DeltaXcyc}) is just for one cross-section. To acquire $\delta \vect{x}_{\text{cyc}}$ for all cross sections, the formula governing its evolution along a cycle in cylindrical coordinates is concluded as below~(\ref{eq:deltaXcyc_evolve_cylindrical})
\begin{subequations}
\begin{align}
&\cancel{
\frac{\rmd}{\rmd t} \delta_{\perp} \vect{x}_{\text{cyc}} 
 = \nabla \vect{B} \cdot \delta_{\perp} \vect{x}_{\text{cyc}} 
+  \delta \vect{B} },
\label{eq:deltaXcyc_evolve_Cartesian} \\
& \frac{\rmd}{\rmd \phi} \delta \vect{x}_{\text{cyc}} 
= \frac{\partial (R\Bpol / B_\phi)}{\partial (R,Z)}\cdot \delta \vect{x}_{\text{cyc}} 
+ \delta \frac{ R\vect{B}_{\text{pol}} }{ B_\phi }.
\label{eq:deltaXcyc_evolve_cylindrical}
\end{align}        
\end{subequations}
\begin{proof} To derive Eq.~\eqref{eq:deltaXcyc_evolve_cylindrical}, differentiate both sides of Eq.~\eqref{eq:DeltaXcyc} \textit{w.r.t.} $\phi$, then
\begin{align*}
    \frac{\rmd }{\rmd \phi} \delta\vect{x}_{\text{cyc}} 
&= - \underbrace{
    \frac{\rmd }{\rmd \phi} \left[
    \calDPm - \matr{I}
    \right]^{-1}
}_{\mathclap{
    \text{utilize a fact for matrix-valued function derivative: } (\matr{A}^{-1})^\prime = \matr{A}^{-1} (-\matr{A}^{\prime}) \matr{A}^{-1}
}} \cdot \delta \mathcal{P}^m
    - \left[
        \calDPm - \matr{I}
    \right]^{-1} \cdot \frac{\rmd }{\rmd \phi} \delta \mathcal{P}^m
\nonumber \\
    &= - \left[ \calDPm - \matr{I} \right]^{-1} (-\frac{\rmd }{\rmd \phi} \calDPm) \left[ \calDPm - \matr{I} \right]^{-1}\cdot \delta \mathcal{P}^m
\nonumber \\
    &\quad - \left[ \calDPm - \matr{I} \right]^{-1} \cdot \frac{\rmd }{\rmd \phi} \delta \mathcal{P}^m
\nonumber \\
    &= - \left[ \calDPm - \matr{I} \right]^{-1} 
\nonumber \\
    &\qquad\cdot\left(
        (-\frac{\rmd }{\rmd \phi} \calDPm) \left[ \calDPm - \matr{I} \right]^{-1}\cdot \delta \mathcal{P}^m
        + \frac{\rmd }{\rmd \phi} \delta \mathcal{P}^m
    \right)
\nonumber
\intertext{(Note that $\frac{\rmd }{\rmd \phi} \calDPm = \matr{A}\cdot \calDPm - \calDPm\cdot \matr{A}$ by Eq.~\eqref{eq:DPm_evolution} and the expression for $\frac{\rmd }{\rmd \phi} \delta\calPm $ can be found in Eq.~\eqref{eq:deltaPm_evolve_cylindrical})}
&= - \matr{A} \left[ \calDPm - \matr{I} \right]^{-1} \cdot \delta \mathcal{P}^m + \delta \frac{R\Bpol }{B_\phi}
\nonumber \\
&=   \matr{A} \delta \xcyc + \delta \frac{R\Bpol }{B_\phi}
\end{align*}
\end{proof}
A concise expression like Eq.~\eqref{eq:deltaXcyc_evolve_Cartesian} for $\frac{\rmd}{\rmd \phi} ~ \delta_{\perp} \vect{x}_{\text{cyc}}$ in an $N$-D flow as a counterpart of Eq.~(\ref{eq:deltaXcyc_evolve_cylindrical}) has not been found. What the authors deduced is too verbose to be a formula, of which the main complexity comes from the fact that $\hatb$ changes along the cycle so the equation involves $\vect{B}\cdot\nabla\hatb$ and becomes difficult to reduce, explained in great details in Appendix \ref{ap:ddt_delta_xcyc_irreducible_expression}. If the requirement of $\delta \xcyc$ being always perpendicular to $\hatb$ during evolution can be dropped, one can still evolve the variation along the cycle by reusing Eq.~\eqref{eq:progress_deltaX_cartesian}
\begin{align}
&
\frac{\rmd}{\rmd t} \delta \vect{x}_{\text{cyc}} 
 = \nabla\vect{B} \cdot \delta \vect{x}_{\text{cyc}} 
+  \delta \vect{B} ,
\label{eq:deltaXcyc_evolve_Cartesian_not_necessarily_perp}
\end{align}
setting the initial condition $\delta\xcyc(t=0)$ to equal the perpendicular shift, $\delta_\perp \xcyc (t=0)$.

Hereafter, the full-period Jacobian on a cycle and its change under perturbation will be examined. As a first step,  with the whole system $\calB$ considered as an argument, the equation $\partial_{\phi_e} \calDXpol (\vect{x}_{0,\text{pol}}, \phi_s, \phi_e ) 
= \matr{A} ~ \calDXpol $~(\ref{eq:DXpol_timeforward}) is complicated into 
\begin{align}
& \frac{\partial }{\partial \phi_e} \calDXpol [\calB ] (\vect{x}_{\text{cyc}} [ \calB  ] (\phi_s), \phi_s, \phi_e )  
\\
&= \matr{A} [\calB] ( \vect{x}_{\text{cyc}} [\calB] (\phi_e) ,\phi_e ) 
\cdot \calDXpol [\calB] ( \vect{x}_{\text{cyc}}[\calB] (\phi_s), \phi_s, \phi_e ) ,
\nonumber
\end{align}
of which an integration in $\phi_e$ after imposed $\Delta\calB\cdot \frac{\rmd}{\rmd \calB}$ gives $(\Delta \calB \cdot \frac{\rmd}{\rmd \calB} ) \calDXpol$ on LHS. Note that this is a total derivative considering the influence from the $\vect{x}_{\text{cyc}}$ shift. 
The $\calDPm$ evolution formula $ \frac{\rmd}{\rmd \phi} ~\calDPm  = \Bigl[ \matr{A} , \calDPm \Bigr] 
$ \eqref{eq:DX_timeforward} can be similarly processed, which yields $ (\Delta \calB \cdot \frac{\rmd}{\rmd \calB} ) \calDPm$ on all sections. The relevant deduction is put in Appendix \ref{ap:change_of_full_period_Jacobian}.


By the \emph{$\calDPm$ evolution formula}, the $\calDPm$ of X-cycles at the edge island chain in the standard configuration of Wendelstein 7-X on all $\phi$-sections have been calculated and shown in Fig.~\ref{fig:W7X_Poincare_3D} by the eigenvectors. Let the vacuum field of the upper control coils be the perturbation, \textit{i.e.} $\calB + \Delta\calB = \calB_{\text{std}} + k\delta\calB_{\text{pert}}$. 
The computed $(\Delta \calB \cdot \frac{\rmd}{\rmd \calB} ) \calDPm$ are verified by results acquired in a brute-force way (that is, simply relocate the X-cycle and recalculate $\calDPm$), as shown in Fig.~\ref{fig:DPm_change_Xcyc_of_W7X}(c), where the eigenvectors of $\calDPm$ at different $k$ are drawn as ribbons to present the change tendency vividly. The yellow dashed edges corresponding to the computed $(1+ \Delta \calB \cdot \frac{\rmd}{\rmd \calB} ) \calDPm$ match the black edge lines corresponding to $k=1$ acquired in the brute-force way as expected. 

The theory can be applied not only to fix an X-cycle at a specified position in the divertor region, but also to maintain the eigenvalues of $\calDPm$ close to unity during operation. This ensures longer field line connection lengths, which facilitates detachment, thereby counteracting potentially detrimental effects induced by the plasma response field~\cite{hudson2002, loizu2017, zhou2022, knieps2022plasma}.

\section{Conclusion and discussion}





A \textit{functional perturbation theory}~(FPT) on the shifts of orbits/trajectories and periodic orbits/cycles has been developed, providing a deeper understanding of the nature of chaotic fields. This understanding allows fusion machines to move beyond relying solely on plasma diffusion to increase the power decay length $\lambda_{q}$ of the scrape-off layer~\cite{eich2011}. Instead, a proactive stimulation of chaotic fields can be employed to disperse heat flux in midair before it reaches the target plate. Moreover, this theory has implications beyond fusion research, offering insights into the behaviour and sensitivity to perturbation of complex systems in various domains. Agile and accurate prediction of orbit and periodic orbit shifts under perturbation can guide practitioners in quickly identifying the desired perturbation to the system.


\begin{acknowledgments}
This work has been carried out within the framework of the EUROfusion Consortium, funded by the European Union via the Euratom Research and Training Programme (Grant Agreement No 101052200 –EUROfusion). Views and opinions expressed are however those of the author(s) only and do not necessarily reflect those of the European Union or the European Commission. Neither the European Union nor the European Commission can be held responsible for them. 
\end{acknowledgments}

\appendix

\section{High-order shifts of orbits and periodic orbits}
\label{ap:high_order_shifts}

\begin{figure}
\includegraphics[width=1.05\linewidth]{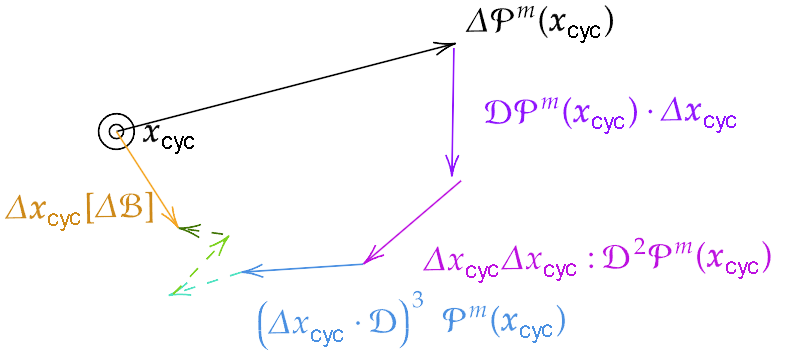}
\caption{\label{fig:xcyc_shift_considering_high_order} Cartoon to show how to calculate the shift of a periodic point for a map under perturbation}
\end{figure}

In the main text, a simple deduction for $\delta\xcyc$ is presented with only $\delta\calPm$ and $\calDPm$ involved. With higher order effects taken into account, the equation (12) is complicated into Eq.~\eqref{eq:match_equation_with_high_orders} as shown below and in Fig.~\ref{fig:xcyc_shift_considering_high_order},
\begin{align}
\Delta \vect{x}_\text{cyc} 
&= \Delta\mathcal{P}^m (\vect{x}_\text{cyc}) 
+ \calDPm(\vect{x}_\text{cyc})\cdot \Delta \vect{x}_\text{cyc} ,
\tag{12 revisited}\\
\Delta \vect{x}_\text{cyc} 
&= \Delta\mathcal{P}^m (\xcyc) 
+ \left(
\Delta\xcyc \cdot \mathcal{D}
\right) \mathcal{P}^m(\xcyc)
\nonumber \\ &
\quad + \frac{1}{2!} \left(
\Delta\xcyc\Delta\xcyc \cddot \mathcal{D}^2
\right) \calPm(\xcyc)
\dots
\nonumber\\
&= \Delta\mathcal{P}^m (\xcyc) 
+ \sum_{k=1}^{\infty} \left(
\Delta\xcyc^k \cdddot_{(k)} \mathcal{D}^k
\right) \mathcal{P}^m(\xcyc) ,
\label{eq:match_equation_with_high_orders}
\end{align}
where $\Delta\xcyc^k \cdddot_{(k)} \mathcal{D}^k := (\Delta\xcyc \cdot \mathcal{D})^k$. Expand $\Delta\xcyc$ and $\Delta \calPm $ by their variations of different orders:  
\begin{align*}
\Delta\xcyc =& \delta \xcyc + \frac{1}{2!} \delta^2 \xcyc + \frac{1}{3!} \delta^3 \xcyc + \cdots ,
\qquad 
\nonumber \\ 
\Delta \calPm = &
\delta \calPm + \frac{1}{2!} \delta^2 \calPm 
+ \frac{1}{3!}\delta^3 \calPm + \cdots 
\end{align*}
Similarly, $\mathcal{D}^k\calPm$ also needs to be expanded:
\begin{align*}
\calD^k\calPm [\calP+\Delta \calP] (\xcyc[\calP]) = \calD^k\calPm [\calP](\xcyc)  
+ \delta \calD^k \calPm[\calP; \Delta\calP](\xcyc)
\nonumber \\
+ \frac{1}{2!} \delta^2 \calD^k \calPm[\calP; \Delta\calP, \Delta\calP](\xcyc)
+ \cdots
\end{align*}

Then, to determine $\delta^k \xcyc$, collect terms in Eq.~\eqref{eq:match_equation_with_high_orders} that contain an equal number of variations, denoted by $k$, to construct the equation of $k^{\text{th}}$ order. Each functional variation $\delta^{k^\prime}$ of order $k^\prime$ contributes $k^\prime$ to this count. For example, the first four orders equations are
\begin{align}
    &\delta\xcyc = \delta\calPm + \delta\xcyc \cdot \calDPm 
,  \label{eq:delta2_xcyc_collect} \\
    &\frac{1}{2!} \delta^2 \xcyc 
    = \frac{1}{2!}\delta^2 \calPm 
    + \frac{1}{2!} \delta^2 \xcyc \cdot \calDPm
    + \delta\xcyc \cdot \delta \calDPm
\nonumber \\
    &\qquad + \frac{1}{2!} \delta\xcyc\delta\xcyc\cddot \mathcal{D}^2\calPm
    , \nonumber
\\
    &\frac{1}{3!} \delta^3 \xcyc 
    = \frac{1}{3!}\delta^3 \calPm 
    + \frac{1}{3!} \delta^3 \xcyc \cdot \calDPm
    + \frac{1}{2!} \delta^2 \xcyc \cdot \frac{1}{1!} \delta  \calDPm
\nonumber \\
    &\qquad
    + \frac{1}{1!} \delta   \xcyc \cdot \frac{1}{2!} \delta^2\calDPm 
    \label{eq:delta3_xcyc_collect} 
    + \frac{1}{2!} \frac{1\times 2}{2!} \delta^2 \xcyc \delta \xcyc \cddot \calD^2 \calPm
\\
    & \qquad + \frac{1}{2!} \delta\xcyc \delta\xcyc \cddot \delta \calD^2\calPm
    + \frac{1}{3!} \delta\xcyc\delta\xcyc\delta\xcyc\cdddot \mathcal{D}^3\calPm , 
\nonumber \\
    &\frac{1}{4!} \delta^4 \xcyc 
    = \frac{1}{4!}\delta^4 \calPm 
    + \frac{1}{4!} \delta^4 \xcyc \cdot \calDPm
    + \frac{1}{3!} \delta^3 \xcyc \cdot \delta \calDPm
\nonumber \\
    &\qquad + \frac{1}{2!} \delta^2 \xcyc \cdot \frac{1}{2!} \delta^2 \calDPm
    + \frac{1}{1!} \delta   \xcyc \cdot \frac{1}{3!} \delta^3 \calDPm 
\nonumber \\
    &\qquad + \frac{1}{2!} \frac{1\times 2}{3!} \delta^3 \xcyc \delta \xcyc \cddot \calD^2 \calPm
    + \frac{1}{2!} \frac{1}{2!} \delta^2 \xcyc \frac{1}{2!} \delta^2 \xcyc \cddot \calD^2 \calPm
\nonumber \\
    &\qquad + \frac{1}{2!} \frac{2}{2!} \delta^2 \xcyc \delta\xcyc \cddot \delta \calD^2\calPm
    + \frac{1}{2!}  \delta \xcyc \delta\xcyc \cddot \frac{1}{2!} \delta^2 \calD^2\calPm
\nonumber \\
    &\qquad+ \frac{1}{3!} \frac{3}{2!} \delta^2\xcyc\delta\xcyc\delta\xcyc\cdddot \mathcal{D}^3\calPm
    + \frac{1}{3!} \delta\xcyc\delta\xcyc\delta\xcyc\cdddot \delta \mathcal{D}^3\calPm
\nonumber \\
    &\qquad+ \frac{1}{4!} \delta\xcyc\delta\xcyc\delta\xcyc\delta\xcyc\cdddot_{(4)} \mathcal{D}^4\calPm , 
\label{eq:delta4_xcyc_collect} 
\end{align}
From the pattern of terms appearing in the above equations, one can conclude the $n^\text{th}$ order formula for the shift of a cycle as below,
\begin{widetext}
\begin{align}
    \frac{1}{n!} \delta^n \xcyc 
& = 
    \sum_{\substack{
        \{p_i\}, \{n_i\}, n_{\calP} ~ \text{such that}\\ 
        n_1 p_1 + \cdots + n_d p_d + n_{\calP} =n
    } } 
    \binom{p^+}{p_1, \dots, p_d}
    (\frac{1}{n_1 !} \delta^{n_1} \xcyc)^{p_1}\cdots (\frac{1}{n_d !} \delta^{n_d} \xcyc)^{p_d} \cdddot_{(p^+)} 
    \frac{
    \delta^{n_{\calP}} \calD^{p^+} \calPm 
    }{n_{\calP}! p^+!} 
    \label{eq:deltan_cyc_collect}
\end{align}    
\end{widetext}
where 
\begin{align*}
    & p_i \geq 1,\qquad n_1 > n_2 > \cdots > n_d \geq 1,
\\
    & p^+ := p_1 + p_2 + \cdots + p_d, \quad\text{ $d$ is the number of powers,}
\\
    & n \geq p^+ \geq 0, \quad \text{ and } n_{\calP} \geq 0  .
\end{align*}
Note that on the RHS of Eq.~\eqref{eq:deltan_cyc_collect} there exists a term $\delta^n\xcyc\cdot \calDPm / n!$. To solve for $\delta^n \xcyc$, one needs to move this term to LHS and multiply both sides by a matrix $-(\calDPm - \matr{I})^{-1}$. 

One can also directly impose $\Delta\calB\cdot \rmd / \rmd \calB$ on both sides of the equation of order $k$ to acquire the equation of order $k+1$ recursively. For example, the equation of second order can be deduced from that of first order 
\begin{align}
\delta\xcyc &= - (\calDPm - \matr{I})^{-1} \delta \calPm.
\tag{\ref{eq:DeltaXcyc} revisited}
\end{align}
Note that for a matrix-valued function $\matr{A}(x)$,  $(\matr{A}^{-1})^\prime = - \matr{A}^{-1} \matr{A}^\prime \matr{A}^{-1} $,
\begin{align}
\delta^2 \xcyc &= 
(\calDPm - \matr{I})^{-1} ~
\underbrace{
    (\Delta\calB\cdot \frac{\rmd}{\rmd \calB} )\calDPm
}_{= \delta\xcyc\cdot\mathcal{D}^2\calPm+\delta\calDPm}~
\underbrace{
    (\calDPm - \matr{I})^{-1}
    \delta \calPm
    }_{=-\delta \xcyc}
\nonumber \\ & -
(\calDPm - \matr{I})^{-1}~ 
\underbrace{
    (\Delta\calB\cdot \frac{\rmd}{\rmd \calB} )\delta\calPm
}_{
    \delta^2\calPm + \delta\xcyc\cdot \delta\calDPm
}. \label{eq:delta2_xcyc_direct_diff}
\end{align}
Since $\calPm$ here is evaluated at the periodic point $\xcyc$ (obviously dependent on $\calB$), the directional functional derivative  $\Delta \calB \cdot \rmd / \rmd \calB$ imposed on $\delta^{k_\delta}\mathcal{D}^{k_D} \calPm$ not only brings a term $\delta^{k_\delta +1}\mathcal{D}^{k_D} \calPm$ but also $\delta \xcyc \cdot \delta^{k_\delta }\mathcal{D}^{k_D+1} \calPm$. Equations \eqref{eq:delta2_xcyc_collect} and \eqref{eq:delta2_xcyc_direct_diff} lead to the same result:
\begin{align}
\delta^2 \xcyc =
- (\calDPm - \matr{I})^{-1}
\left(
\delta^2\calPm
+ 2 \delta\xcyc\cdot \delta\calDPm
+\delta\xcyc\delta\xcyc\cddot \calD^2\calPm
\right),
\end{align}
which can be further imposed $\Delta\calB\cdot \rmd / \rmd \calB$ on both sides to acquire expressions for the third and fourth order,
\begin{widetext}
\begin{align}
\delta^3 \xcyc = &
-(\calDPm - \matr{I})^{-1} -(\Delta\calB\cdot \frac{\rmd}{\rmd \calB} ) (\calDPm - \matr{I})
\overbrace{
(\calDPm - \matr{I})^{-1}
\left(
\delta^2\calPm
+ 2 \delta\xcyc\cdot \delta\calDPm
+\delta\xcyc\delta\xcyc\cddot \calD^2\calPm
\right)
}^{-\delta^2 \xcyc} \\
& -(\calDPm - \matr{I})^{-1}
\Big(
\delta \xcyc \cdot \delta^2 \calD \calPm
+ \delta^3 \calPm
+ 2 \delta^2 \xcyc \cdot \delta \calD\calPm
+ 2 \delta \xcyc \delta \xcyc \cddot \delta \calD^2 \calPm 
+ 2 \delta \xcyc \cdot \delta^2 \calD \calPm
\nonumber \\
& \qquad\qquad
+ 2 \delta^2 \xcyc \delta \xcyc \cddot \calD^2 \calPm
+ \delta\xcyc\delta\xcyc\delta\xcyc\cdddot \calD^3 \calPm
+ \delta\xcyc\delta\xcyc\cddot \delta\calD^2 \calPm
\Big) 
\nonumber \\
= & -(\calDPm - \matr{I})^{-1}
\Big(
\delta^3 \calPm
+ 3 \delta^2 \xcyc \cdot \delta \calD\calPm
+ 3 \delta \xcyc \delta \xcyc \cddot \delta \calD^2 \calPm 
\\
& \qquad\qquad
+ 3 \delta \xcyc \cdot \delta^2 \calD \calPm
+ 3 \delta^2 \xcyc \delta \xcyc \cddot \calD^2 \calPm
+ \delta\xcyc\delta\xcyc\delta\xcyc\cdddot \calD^3 \calPm
\Big) \nonumber
\end{align}
\begin{align}
\delta^4 \xcyc 
= & -(\calDPm - \matr{I})^{-1}
\Big(
\delta^4 \calPm
+ 4 \delta^3 \xcyc \cdot \delta \calD\calPm
+ 6 \delta^2 \xcyc \cdot \delta^2\calD\calPm
+ 6 \delta \xcyc \delta \xcyc \cddot \delta^2 \calD^2 \calPm 
\\
& \qquad\qquad
+ 4 \delta \xcyc \cdot \delta^3 \calD \calPm
+ 4 \delta^3 \xcyc \delta \xcyc \cddot \calD^2 \calPm
+ 6 \delta^2 \xcyc\delta\xcyc\delta\xcyc\cdddot \calD^3 \calPm
+ 12 \delta^2 \xcyc\delta\xcyc\cddot \delta\calD^2 \calPm 
\nonumber \\
& \qquad\qquad
+ 3 \delta^2 \xcyc\delta^2 \xcyc \cddot \calD^2\calPm
+ 4 \delta\xcyc \delta\xcyc \delta\xcyc \cdddot \delta \calD^3\calPm
+ \delta\xcyc \delta\xcyc \delta\xcyc\delta\xcyc \cdddot_{(4)} \calD^4\calPm
\Big) 
\nonumber
\end{align}
\end{widetext}

Similarly, the progression of trajectory variations (general trajectories meant here, which do not need to be periodic) such as those of the first three orders are described by Eq.~(\ref{eq:progress_deltaX_cartesian}, \ref{eq:delta2X_timeforward}, \ref{eq:delta3X_timeforward}). One can obtain the $n^\text{th}$ order trajectory variation progression formula as below in a manner similar to Eq.~\eqref{eq:deltan_cyc_collect}, 
\begin{widetext}
\begin{align}
\frac{1}{n!}
\frac{\partial }{\partial t} \delta^n \vect{X} 
& = \sum_{\substack{
\{n_i\}, \{p_i\}, n_{\vect{B}} ~\text{such that} \\
n_1 p_1 + \cdots + n_d p_d + n_{\vect{B}} =n 
} } 
\binom{p^+}{p_1, \dots, p_d}
(\frac{1}{n_1 !} \delta^{n_1} \vect{X})^{p_1}\cdots (\frac{1}{n_d !} \delta^{n_d} \vect{X} )^{p_d} \cdddot_{(p^+)} \frac{1}{ p^+ ! n_{\vect{B}}! } \nabla^{p^+} \delta^{n_{\vect{B}}}  \vect{B} ,
\label{eq:nth_variation_progerssion_formula}
\end{align}
\end{widetext}
where 
\begin{align*}
    & p_i \geq 1,\qquad n_1 > n_2 > \cdots > n_d \geq 1,
\\
    & p^+ := p_1 + p_2 + \cdots + p_d, \quad\text{ $d$ is the number of powers,}
\\
    & n \geq p^+ \geq 0, \quad \text{ and } n_{\vect{B}} \geq 0  .
\end{align*}
The formula~\eqref{eq:nth_variation_progerssion_formula} can be easily migrated from $N$-D continuous-time flows to $N$-D discrete-time maps, \textit{e.g.} $\calP$, to describe the progression of discrete orbit variations by replacing symbols as below
\begin{align*}
\frac{\partial}{\partial t }
 \vect{X} (\vectx_0, t) \text{  on LHS  } &\mapsto \vect{X}(\vect{x}_{0}, k+1) = \calP^{k+1} (\vectx_0),
\\
\vect{X} (\vectx_0, t) \text{ on RHS} &\mapsto \vect{X}(\vectx_{0}, k) = \calP^{k} (\vectx_0),,
\\
\vect{B} & \mapsto \calP,   
\\ \nabla\vect{B} &\mapsto \calD\calP.
\end{align*}

\section{Evolution of $\delta_\perp \xcyc$ along a cycle}
\label{ap:ddt_delta_xcyc_irreducible_expression}
As stated in the main text, the evolution of the perpendicular shift of a cycle, $\delta_\perp \xcyc$, along a cycle is too verbose to be called a formula. The deduction details are shown below
\begin{align}
\delta_{\perp} \vect{x}_{\text{cyc}}  
= &
\left[
    \matr{I} - (\matr{I} - \hat{\vect{b}}\hat{\vect{b}}^T) \cdot \calDXT
\right]^{-1} \cdot
(\matr{I}-\hat{\vect{b}}\hat{\vect{b}}^T)\cdot \delta \vect{X}_T 
\tag{15 revisited} 
\end{align}
\begin{align}
\frac{\rmd}{\rmd t} \delta_{\perp} \vect{x}_{\text{cyc}}  
= &
\frac{\rmd}{\rmd t}
\left[
    \matr{I} - (\matr{I} - \hat{\vect{b}}\hat{\vect{b}}^T) \cdot \calDXT
\right]^{-1} \cdot
(\matr{I}-\hat{\vect{b}}\hat{\vect{b}}^T)\cdot \delta \vect{X}_T  
\nonumber \\
& +
\left[
    \matr{I} - (\matr{I} - \hat{\vect{b}}\hat{\vect{b}}^T) \cdot \calDXT
\right]^{-1} \cdot
\frac{\rmd}{\rmd t} (\matr{I}-\hat{\vect{b}}\hat{\vect{b}}^T)\cdot \delta \vect{X}_T  
\nonumber\\
& +
\left[
    \matr{I} - (\matr{I} - \hat{\vect{b}}\hat{\vect{b}}^T) \cdot \calDXT
\right]^{-1} \cdot
(\matr{I}-\hat{\vect{b}}\hat{\vect{b}}^T)\cdot \frac{\rmd}{\rmd t}\delta \vect{X}_T  \\
= &
-\left[
    \matr{I} - (\matr{I} - \hat{\vect{b}}\hat{\vect{b}}^T) \cdot \calDXT
\right]^{-1} \cdot
\frac{\rmd}{\rmd t} 
\left[
    \matr{I} - (\matr{I} - \hat{\vect{b}}\hat{\vect{b}}^T) \cdot \calDXT
\right] \nonumber \\
& \qquad
\cdot\underbrace{
\left[
    \matr{I} - (\matr{I} - \hat{\vect{b}}\hat{\vect{b}}^T) \cdot \calDXT
\right]^{-1} \cdot
(\matr{I}-\hat{\vect{b}}\hat{\vect{b}}^T)\cdot \delta \vect{X}_T  
}_{ = \delta_\perp \xcyc }
\nonumber \\
& +
\left[
    \matr{I} - (\matr{I} - \hat{\vect{b}}\hat{\vect{b}}^T) \cdot \calDXT
\right]^{-1} \cdot
\frac{\rmd}{\rmd t} (\matr{I}-\hat{\vect{b}}\hat{\vect{b}}^T)\cdot \delta \vect{X}_T  
\nonumber\\
& +
\left[
    \matr{I} - (\matr{I} - \hat{\vect{b}}\hat{\vect{b}}^T) \cdot \calDXT
\right]^{-1} \cdot
(\matr{I}-\hat{\vect{b}}\hat{\vect{b}}^T)\cdot \frac{\rmd}{\rmd t}\delta \vect{X}_T  
\end{align}
where $\dfrac{\rmd}{\rmd t} (\matr{I} - \hat{\vect{b}}\hat{\vect{b}}^T) 
= - \vect{B}\cdot \nabla (\hat{\vect{b}}\hat{\vect{b}}^T)
= - B (\vect{\kappa}\hat{\vect{b}}^T + \hat{\vect{b} } \hat{\vect{\kappa}}^T) $, while $\dfrac{\rmd}{\rmd t} \delta \vect{X}_T = \matr{A} ~\delta \vect{X}_T - (\calDXT - \matr{I}) \delta \vect{B}$, which is Eq.~\eqref{eq:deltaXT_evolve_Cartesian} in the main text. Even though these two derivatives have clear expressions, the RHS seems hard to reduce further.

\section{\\The change of full-period Jacobian $(\Delta \calB \cdot \rmd / \rmd \calB )~ \calD\calPm$}
\label{ap:change_of_full_period_Jacobian}
For the full-period Jacobian of a cycle and its change under perturbation discussed in the main text, the derivation and detailed formulae are shown here. As a first step,  with the whole system $\calB$ considered as an argument, the $\calD\vect{X}_\text{pol}$ progression equation $\partial_{\phi_e} \calDXpol (\vect{x}_{0,\text{pol}}, \phi_s, \phi_e ) 
= \matr{A} ~ \calDXpol $, Eq.~\eqref{eq:DXpol_timeforward} in the main text, is complicated into 
\begin{align}
& \frac{\partial }{\partial \phi_e} \calDXpol [\calB ] (\vect{x}_{\text{cyc}} [ \calB  ] (\phi_s), \phi_s, \phi_e )  \label{eq:DPpol_timeforward_variantB} \\
&= \matr{A} [\calB] ( \vect{x}_{\text{cyc}} [\calB] (\phi_e) ,\phi_e ) 
~ \calDXpol [\calB] ( \vect{x}_{\text{cyc}}[\calB] (\phi_s), \phi_s, \phi_e ) ,
\nonumber
\end{align}
which after imposed $\Delta\calB\cdot \frac{\rmd}{\rmd \calB} \calB$ converts to
\begin{widetext}
\begin{align}
\frac{\partial}{\partial \phi_e}
& 
(\Delta \calB \cdot \frac{\rmd}{\rmd \calB} ) ~
\calDXpol 
[\calB] 
( \vect{x}_{ \tikzmark{DXpol_D_} \text{cyc}}  [\calB] (\phi_s), \phi_s, \phi_e ) 
= 
\nonumber \\
\matr{A} 
&
\left(
    \Delta \calB \cdot \frac{ \rmd}{\rmd \calB}
\right)~
\calDXpol 
 +
\left( 
   \Delta \calB \cdot \frac{\tikzmark{A_deltaB} \delta}{\delta \calB} 
   + \left( 
        \Delta \calB \cdot \frac{ \delta \vect{x}_{\text{cyc}} [\mathcal{B}](\phi_e) }{  \delta \mathcal{B}   }
    \right) \cdot \frac{\partial  }{ \partial (R \tikzmark{A_partialRZ} ,Z)} 
\right) ~ 
\Bigl(
    \matr{A} 
    [\tikzmark{A_deltaB_} \calB] 
    ( \vect{x}_{\tikzmark{A_partialRZ_} \text{cyc}} [\calB] (\phi_e) ,\phi_e ) 
\Bigr)~ \calDXpol 
\tikz[overlay,remember picture]{
   \draw[->,up square arrow] (A_deltaB.north) to (A_deltaB_.north);
   \draw[->,down square arrow] (A_partialRZ.south) to (A_partialRZ_.south);
}
\label{eq:deltaDXpol_timeforward}
\end{align}
\end{widetext}
giving $(\Delta \calB \cdot \frac{\rmd}{\rmd \calB} ) \calDXpol$ on LHS after an integration in $\phi_e$. Note that this is a total derivative considering $\vect{x}_{\text{cyc}}$ shift. 
Furthermore, to avoid repeated computation for $(\Delta \calB \cdot \frac{\rmd}{\rmd \calB} ) \calDXpol$ at all $\phi$-sections, one can similarly process the $\calDPm$ evolution formula 
$ \frac{\rmd}{\rmd \phi}~ \calDPm  = \Bigl[ \matr{A} , \calDPm \Bigr] $
 to let the $\calB$ argument be explicit, 
 
\begin{align}
& 
\frac{\rmd }{\rmd \phi} \calDPm [\calB ] (\vect{x}_{\text{cyc}} [ \calB  ] (\phi), \phi ) 
\label{eq:DPm_evolve_variantB} \\
&= 
\Bigl[
 \matr{A} [\calB] ( \vect{x}_{\text{cyc}} [\calB] (\phi) ,\phi ), ~~
 \calDPm [\calB] ( \vect{x}_{\text{cyc}}[\calB] (\phi), \phi) 
\Bigr],  
\nonumber
\end{align}

which after imposed $\Delta\calB\cdot \rmd / \rmd \calB$ on both sides becomes, 
\begin{widetext}
\begin{align}
\frac{\rmd }{\rmd \phi} 
& (\Delta \calB \cdot \frac{\rmd}{\rmd \calB} ) 
~\calDPm [\calB ] (\vect{x}_{\text{cyc}} [ \calB  ] (\phi), \phi )
= \nonumber \\
\Bigl[
     \matr{A}, ~~
     &
    \left(
    \Delta \calB \cdot \frac{ \rmd}{\rmd \calB}
\right) \calDPm  
\Bigr]
+ 
\Bigl[
    \left( 
       \Delta \calB \cdot \frac{ \tikzmark{A_deltaB2} \delta}{  \delta \calB} 
       + \left( 
            \Delta \calB \cdot \frac{ \delta \vect{x}_{\text{cyc}} [\mathcal{B}](\phi) }{  \delta \mathcal{B}   }
        \right) \cdot \frac{\partial  }{\partial (R  
         \tikzmark{A_partialRZ2} ,Z)} 
    \right)
\matr{A} [ \tikzmark{A_deltaB2_} \calB] ( 
\vect{x}_{ 
\tikzmark{A_partialRZ2_} \text{cyc}
} [\calB] (\phi) ,\phi ), ~~  \calDPm 
\Bigr] ,
\tikz[overlay,remember picture]{
   \draw[->,up square arrow] (A_deltaB2.north) to (A_deltaB2_.north);
   \draw[->,down square arrow] (A_partialRZ2.south) to (A_partialRZ2_.south);
}
\label{eq:delta_DPm_evolve_variantB}
\end{align} 
\end{widetext} of which an integration in $\phi$ yields $ (\Delta \calB \cdot \frac{\rmd}{\rmd \calB} ) \calDPm$ on all sections. The initial condition of 
this ODE system is $(\Delta \calB \cdot \frac{\rmd}{\rmd \calB} ) \calDXpol(\phi_s, \phi_e=\phi_s + 2m\pi)$, coming from a $2m\pi$ integration in $\phi_e$ of Eq.~(\ref{eq:deltaDXpol_timeforward}).

Only the second term on the RHS of Eq.~(\ref{eq:deltaDXpol_timeforward}), that is $(\Delta \calB \cdot \delta / \delta \calB )~ \matr{A} $, may need to be explained to readers. The other terms are comparatively easy to handle.
Components of it are shown in a matrix form
\begin{align}
& (\Delta \calB \cdot \frac{\delta}{\delta \calB} )~ \matr{A} \left[\calB \right]  (\vect{x}_{\text{cyc}} \left[\calB \right](\phi_s) , \phi_s, \phi_e ) 
\intertext{(Note that $\delta / \delta \calB$ is a partial derivative, so $\xcyc$ does not matter here.)}
&= 
(\Delta \calB \cdot \frac{\delta}{\delta \calB} )~ 
    \begin{bmatrix}
    \frac{B_{R}}{B_\phi} + R \left( \frac{ \partial_R B_R}{B_\phi}  - \frac{B_R \partial_R B_\phi}{B_\phi^2}\right) & 
    R \left( \frac{ \partial_Z B_R}{B_\phi}  - \frac{B_R \partial_Z B_\phi}{B_\phi^2}\right) 
    \\
    \frac{B_{Z}}{B_\phi} + R \left( \frac{ \partial_R B_Z}{B_\phi}  - \frac{B_Z \partial_R B_\phi}{B_\phi^2}\right) & 
    R \left( \frac{ \partial_Z B_Z}{B_\phi}  - \frac{B_Z \partial_Z B_\phi}{B_\phi^2}\right) 
    \end{bmatrix} 
\nonumber \\
    &= \begin{bmatrix}
    \frac{\delta B_{R}}{B_\phi} - \frac{B_R}{B_\phi^2} \delta B_\phi + (1,1) & (1,2) \\
    \frac{\delta B_{Z}}{B_\phi} - \frac{B_Z}{B_\phi^2} \delta B_\phi + (2,1) & (2,2)
    \end{bmatrix}
\end{align}
\begin{align*}
(1,1) &= 
R \Big( 
    \frac{\partial_R \delta B_R}{B_\phi}
    -\frac{\partial_R B_R }{B_\phi^2} \delta B_\phi 
    -\frac{ \partial_R B_\phi }{B_\phi^2 }\delta B_R 
\\ & \qquad
    - B_R \frac{ 
        (\partial_R \delta B_\phi) B_\phi^2  
        - (\partial_R B_\phi)  2 B_\phi \delta B_\phi
        }{B_\phi^4}\Big) \\ 
(1,2) &= (1,1) \text{ with } \partial_R \text{ replaced by } \partial_Z. \\
(2,1) &= (1,1) \text{ with } B_R \text{ and } \delta B_R \text{ \textit{resp.} replaced by } B_Z \text{ and } \delta B_Z. \\
(2,2) &= (1,1) \text{ with } \partial_R, B_R \text{ and } \delta B_R \text{ \textit{resp.} replaced by } \partial_Z, B_Z \text{ and } \delta B_Z. 
\end{align*}
An alternative expression of $(\Delta \calB \cdot \frac{\delta}{\delta \calB} )~ \matr{A}$ is 
\begin{align}
&
\frac{\partial }{\partial (R,Z) }
\begin{bmatrix}
\delta \dfrac{R\Bpol}{B_\phi}
+
R \dfrac{\partial}{\partial R}
\left(
\delta \dfrac{R\Bpol}{B_\phi}
\right)
&\Big\vert&
R \dfrac{\partial}{\partial Z}
\left(
\delta \dfrac{R\Bpol}{B_\phi}
\right)
\end{bmatrix},
\nonumber
\end{align}
where $\delta \dfrac{R\Bpol}{B_\phi} = (\Delta\calB\cdot \delta / \delta \calB) \dfrac{R\Bpol}{B_\phi} = 
\dfrac{R\delta\Bpol}{B_\phi}
- \dfrac{R\Bpol}{B_\phi^2}\delta B_\phi$.

To migrate equations~(\ref{eq:DPpol_timeforward_variantB}-\ref{eq:delta_DPm_evolve_variantB}) from 3D flows in cylindrical coordinates to $N$-D flows in Cartesian coordinates, replace symbols as below
\begin{align}
\matr{A} \mapsto \nabla\vect{B}, &&        \calDPm \mapsto \calDXT,     &&             \frac{\partial}{\partial (R,Z)} \mapsto \nabla.    
\nonumber
\end{align}
Be aware that the perturbation may change the period of this cycle, so the change in the time integration range may bring new terms in formulae.

To migrate equations~(\ref{eq:DPpol_timeforward_variantB}-\ref{eq:delta_DPm_evolve_variantB}) to $N$-D maps, \textit{e.g.} $\calP$, one needs more than symbol substitution. The counterpart of Eq.~\eqref{eq:DXpol_timeforward} for mapping is 
\begin{align}
    \calD\calP^k (\vectx) 
    = \calD\calP( \calP^{k-1} (\vectx) ) 
    \cdot \calD\calP^{k-1}( \vectx )
    \label{eq:DPk_progress_mapping}
\end{align}
Given an $m$-periodic orbit $\{\xcyc, \calP(\xcyc), \cdots, \calP^{m-1}(\xcyc)\}$, one can compare
\begin{align*}
    \calD\calP^m (\xcyc) 
    =& \calD\calP( \calP^{m-1} (\xcyc) ) 
    \cdot \calD\calP( \calP^{m-2} (\xcyc) )  \cdots
    \calD\calP( \xcyc )
\intertext{and }
    \calD\calP^m ( \calP(\xcyc) ) 
    =& \calD\calP( \calP^{m} (\xcyc) ) 
    \cdot \calD\calP( \calP^{m-1} (\xcyc) )  \cdots
    \calD\calP( \calP(\xcyc) )
\end{align*}
to find the relationship between $\calDPm(\xcyc)$ and $\calDPm( \calP(\xcyc) )$, that is
\begin{align}
    \calD\calPm( \calP(\xcyc) ) 
    & = \calD\calP(\calP^{m}(\xcyc)) \cdot \calDPm(\xcyc)\cdot [ \calD\calP(\xcyc) ]^{-1}
\nonumber \\
    & = \calD\calP(\xcyc) \cdot \calDPm(\xcyc)\cdot [ \calD\calP(\xcyc) ]^{-1},
    \label{eq:DPm_evolution_mapping}
\end{align}
which is the counterpart of $\calDPm$ evolution formula~\eqref{eq:DPm_evolution} for mapping. Upon applied $\Delta\calP\cdot \frac{\rmd}{\rmd \calP}$ on both sides, Eq.~\eqref{eq:DPk_progress_mapping} and \eqref{eq:DPm_evolution_mapping} transform to the progression equation of $(\Delta\calP\cdot\frac{\rmd}{\rmd \calP} ) \calD\calP^k(\xcyc)$,
\begin{align}
    (\Delta\calP\cdot\frac{\rmd}{\rmd \calP} )  \calD\calP^k (\xcyc) 
    = \calD\calP( \calP^{k-1} (\xcyc) ) 
    \cdot (\Delta\calP\cdot\frac{\rmd}{\rmd \calP} ) \calD\calP^{k-1}( \xcyc )
\nonumber \\
    + \underbrace{ (\Delta\calP\cdot\frac{\rmd}{\rmd \calP} ) \left[\calD\calP( \calP^{k-1} (\xcyc) ) \right] }_{\substack{
    = \delta \calD\calP( \calP^{k-1} (\xcyc) ) \\
    + \left[ (\delta\calP^{k-1}+ \delta\xcyc\cdot\calD\calP^{k-1})\cdot\calD\right] \calD\calP( \calP^{k-1} (\xcyc) ) 
    } }
    \cdot \calD\calP^{k-1}( \xcyc ),
\end{align}
and the equation describing the evolution of $(\Delta\calP\cdot\frac{\rmd}{\rmd \calP} ) \calD\calP^k(\xcyc)$ along a discrete periodic orbit,
\begin{align}
    & (\Delta\calP\cdot\frac{\rmd}{\rmd \calP} )  \calD\calPm( \calP(\xcyc) ) 
\nonumber \\
    = & \overbrace{ (\Delta\calP\cdot\frac{\rmd}{\rmd \calP} ) \left[\calD\calP(\xcyc)\right] }^{=\delta\calD\calP+(\delta\xcyc\cdot\calD)\calD\calP} \cdot \calDPm(\xcyc)\cdot [ \calD\calP(\xcyc) ]^{-1}
\nonumber \\
    & + \calD\calP(\xcyc) \cdot (\Delta\calP\cdot\frac{\rmd}{\rmd \calP} ) \calDPm(\xcyc)\cdot [ \calD\calP(\xcyc) ]^{-1}
\nonumber \\
    & + \left[\calD\calP(\xcyc)\right] \cdot \calDPm(\xcyc)\cdot  \underbrace{ (\Delta\calP\cdot\frac{\rmd}{\rmd \calP} )[ \calD\calP(\xcyc) ]^{-1}}_{ \text{utilize } (\matr{A}^{-1})^\prime = \matr{A}^{-1}(-\matr{A}^{\prime})\matr{A}^{-1}  },
\end{align}
which are respectively the counterparts of Eq.~\eqref{eq:deltaDXpol_timeforward} and \eqref{eq:delta_DPm_evolve_variantB}.
\bibliography{main}

\end{document}